# Ancient Indian mathematics needs an honorific place in modern mathematics celebration.

**Steven G. Krantz, Ph.D.**[1]

**Arni S. R. Srinivasa Rao***,[2]

**Abstract:** The Indian tradition in mathematics is long and glorious. It dates back to earliest times, and indeed many of the Indian discoveries from 5000 years ago correspond rather naturally to modern mathematical results.

[1] Professor, Department of Mathematics, Washington University in St. Louis, MO. U.S.A. Email: sk@math.wustl.edu (OR) sgkrantz@gmail.com

[2] Professor and Director, Laboratory for Theory and Mathematical Modeling, Medical College of Georgia, Augusta, GA. U.S.A. Email: arrao@augusta.edu (OR) arni.rao2020@gmail.com
**(\*Corresponding author**)

The government of India announced in 2012 that every year Indian mathematician Srinivasa Ramanujan's birthday of December 22 will be celebrated as national mathematics day. The year 2012 was the great Ramanujan's 125th birth anniversary. The government of India released a commemorative stamp on that occasion as well. Brilliant contributions in number theory by Ramanujan are well known (Andrews and Berndt, 2018 part V - Andrews and Berndt, 2018 part I). Mathematical developments and deeper thinking in India were several thousand years older than Ramanujan. In this comment, we try to recollect some of the ancient Indian mathematics and its mathematicians who did fundamental work in number systems, mathematics of astronomy, calculus, etc.,

In point of fact the Indian tradition in mathematics is long and glorious. It dates back to earliest times, and indeed many of the Indian discoveries from 5000 years ago have analogues in modern mathematical results. It is our contention that the national mathematics day of India should have a very broad scope. It is certainly appropriate to pay homage to the memory of Ramanujan. But let us not forget his many mathematical forebears. And also, the very strong mathematicians, such as Harish Chandra, who have come since his time. What we are really celebrating here is a great intellectual tradition of which all humans can be proud. As the great historians have pointed out, those who forget the past will relinquish the future. Indian Mathematics Day need to represent 5000 years historical milestones of Indian mathematics.

In this comment, we do not discuss the complete development of mathematics during ancient and medieval India. We highlight some of the key connecting points that played a significant role in the world of mathematics much before the works by Ramanujan were revealed.

The origins of the mathematics that emerged in the Indian subcontinent can be seen around the Sulbha sutra period which was around 1200BCE to 500 BCE. During this period the numbers up to $10^{12}$ were counted (in Vedic Sanskrit this number was referred to as *Paradham*). The Vedic period mathematics was confined to the geometry of fire alters and astronomy, and these concepts were used to perform rituals by the priests. Some of the famous names from that era are Budhayana, Apastambha, and Katayana. In the next paragraph, we will describe Sanskrit sounds and their corresponding English numerals. Indian mathematics also introduced the decimal number system that is in use today and the concept of zero as a number. The concept of sine (written as *jaya* in Sanskrit) and cosine (*cojaya*), negative numbers, arithmetic, and algebra were found in ancient Indian mathematics (Dutta and Singh, 1962). The mathematics developed in India was later translated and transmitted to China, East Asia, West Asia, Europe, and Saudi Arabia. The classical period of Indian mathematics was often attributed to the period 200AD to 1400AD during which works of several well-known mathematicians, like Varahamihira, Aryabhata, Brahmagupta, Bhaskara, and Madhava got translated into other languages and transmitted outside the sub-continent.

The number systems present since the Vedic days, especially since the Sukla Yajurveda and their Sanskrit sounds, were as follows: 1 (*Eka*), 10 (*Dasa*), 100 (*Sata*), 1000 (*Sahasra*), $10^4$ (*Aayuta*), $10^5$ (*Laksa or Niyuta*), $10^7$ (*Koti*), $10^{12}$ (*Sanku or Paraardha*), $10^{17}$ (*Maha Sanku*), $10^{22}$ (*Vrnda*), $10^{52}$ (*Samudra*), $10^{62}$ (*Maha-ogha*).

**Table 1.** Sanskrit and English numerals 0 to 9 and their sounds.

(a) Magic square (*Anka-yantra*) for Sun

(b) Magic squares (*Anka-yantra*) for other eight planets

**Figure 1.** Magic Squares for Sun and other eight planets available from ancient Indian literature.

In addition to the number systems of ancient India, still today in India are heard the popular *'Vishnu Sahasra Nama Stotra'* which dates back to the Mahabharata epic. In this, there is a verse that sounds like "*Sahasra Koti Yugadharine Namah*". If we translate this verse, *Sahasra* (as we saw above) means 1000, and *Koti* (as we saw above) means $10^7$, so a simple translation of the phrase *Sahasra Koti* could mean $10^{10}$. The entire phrase has been interpreted in different ways. We do not give our readers all possible interpretations here and confine our commentary only to the number systems.

The depth of work in astronomy, solar systems, geometry, and ground-breaking mathematical calculations by ancient and medieval great scholars in India, for example, Budhayana, Varahamihira, Aryabhatta, Bhaaskara I & II, Pingala, Madhava, Lilavathi, and many more, etc., are well known (see Joseph, 2000, Evans, 2014, and Daji, 1865). It seems that celebrating national mathematics day in India only as part of Ramanujan's birthday is confining the glory and celebration of Indian mathematics to a little over 100 years of the past. Schools and colleges across India have celebrated Ramanujan's birthday for many decades, but that is different from exclusively limiting national day only to the great Ramanujan.

A good deal can be written on ancient scholar's work from India, which can be found, for example in Gupta JC (2019, Bhavana Magazine), Plofker KL (2009), TA Saraswati Amma (1979), Dutta & Singh (1962 and 1938), Bürk, A (1901). In this opinion piece, we highlight only a few of them.

Shulba sutras which were compiled in the Sanskrit language around 2000-800 BCE in India by Budhayana, Manava, Apastamba, Katayana consisted of geometric-shaped fire-alters for performing ancient Indian rituals (Burk, (1901), Plofker (2009)). Some of these sutras also contain the statements of Pythagorean theorems and triples. For example, Apastamba provided the following triples

$$\left\{\begin{array}{c}(3,4,5), (5,12,13), (15,8,17),\\ (12,35,17), (12,35,37), (36, 15, 39)\end{array}\right\}$$

for constructing fire-alters Saraswati Amma (1979).

These sutras can be used to find the approximate value of $\sqrt{2}$ (see Saraswati Amma (1979) and Plofker (2009)), using the expression

$\sqrt{2} \approx 1 + \frac{1}{3} + \frac{1}{3.4} - \frac{1}{3.4.34}$ = 1.41421568627451

In the Vedic period astrology (Jyotisha) of India, the magic squares (Anka-yantras) were used to please and worship nine planets of the solar system (Gupta 2005, Styan 2013). Figure 1 is about ancient magic squares for Sun and the other eight planets in our solar system.

During the 5th century AD, Aryabhata calculated, among many other things, that the moon orbit takes 27.396 days, the value of $\pi = 3.1416$, etc., He is believed to have laid the foundations of the properties of sine and cosine in trigonometry.

Plofker (2009) also writes that the Leibnitz infinite series

$$4\left(1 - \frac{1}{3} + \frac{1}{5} - \frac{1}{7} + \frac{1}{9} - \cdots\right)$$

was known in the works of Indian mathematician Madhava who lived three centuries earlier than Leibnitz.

Bhaskara of 12th Century AD through his famous book *Bijaganita* described the rules of algebraic operations of positive, negative signs, rules of zero (*shuunyah*), and infinity (*antam*). His book also describes obtaining solutions to intermediate equations of the first degree (Ramasubramaniam, Hayashi, and Montelli, 2019). Bhaskara's other book titled *Siddanta siromani* provided a detailed account of and development of Indian astronomy. Computations of the planetary movements, shapes of planets, rotation axis, lunar month days etc., were explained in detail.

What we advocate in this piece is for an exposition of deep-rooted mathematical knowledge in India, and not an exhaustive account of all possible results and conclusions. Several of the ancient texts in the language *Sanskrit* are either lost or preserved in museums.

Srinivasa Ramanujan's work undoubtedly shines as part of modern Indian mathematics but thousands of years ancient mathematical discoveries, the introduction of various branches of pure and applied mathematics needs a proper representation in any celebration of India's contribution to world mathematics.

## Conflicts of interest:

None to be declared.